\documentclass[12pt]{elsarticle}
\usepackage[usenames]{color}
\usepackage{amssymb}
\usepackage{amsmath}
\usepackage{amsthm}
\usepackage{amsfonts}
\usepackage{amscd}
\usepackage{graphicx}
\usepackage[colorlinks=true,
linkcolor=webgreen,
filecolor=webbrown,
citecolor=webgreen]{hyperref}
\definecolor{webgreen}{rgb}{0,.5,0}
\definecolor{webbrown}{rgb}{.6,0,0}
\usepackage{color}
\usepackage{fullpage}
\usepackage{float}
\usepackage{graphics}
\usepackage{latexsym}
\usepackage{epsf}
\usepackage{bookmark}
\usepackage{mathrsfs}
\usepackage{hyperref}
\def\C{\mathbb{C}}
\def\C{\mathbb{C}}
\def\N{\mathbb{N}}

\def\cD{\mathcal D}
\def\cH{\mathcal H}
\def\le{\leqslant}
\def\ge{\geqslant}
\def\hc{\hat c}
\def\pr{^\prime }

\def\eop{\unskip\nobreak\hfil\penalty50\hskip2em\hbox{}\nobreak
\hfill\mbox{$\Box $}\par}
\newtheorem{theorem}{Theorem}[section]

\begin{document}
\begin{frontmatter}
	\title{Discrete orthogonality of the polynomial sequences in the $q$-Askey scheme}

\author{Luis Verde-Star}
 \address{
Department of Mathematics, Universidad Aut\'onoma Metropolitana, Iztapalapa,
 Mexico City,
 Mexico }
\ead{verde@xanum.uam.mx}

	\begin{abstract} 
	We  obtain   weight functions associated with  $q$-linear and $q$-quadratic lattices that yield discrete orthogonality with respect to a quasi-definite moment functional  for  the Askey-Wilson polynomials and all the   polynomial sequences in the q-Askey scheme, with the exception of the continuous $q$-Hermite polynomials.   


\end{abstract}
\end{frontmatter}

\section{Introduction}
The families of basic hypergeometric  orthogonal polynomial sequences included in the q-Askey scheme \cite{Hyp}  have been studied for a long time. Some of those families have a well-known  orthogonality  determined by a weight function $w(x_k)$ defined on a finite or infinite  set of points $x_k$ that belong to a $q$-linear or $q$-quadratic lattice. In such cases the orthogonality of a polynomial sequence $\{u_n(t): n \in \N\}$ is given by 
\begin{equation}\label{discrOrtho}
\sum_{k \in  M} u_n(x_k) u_m(x_k) w(x_k)=  K_n \delta_{n,m}, \qquad  n,m \in M,
\end{equation}
where $M=\{0,1,2,\ldots,N\}$, for some positive integer $N$, or $M=\N$, and  $K_n$ is a  nonzero constant for $n\in M$.

The polynomial families that satisfy certain types of differential or difference equations and have a discrete orthogonality with positive weights $w(x_k)$ have been characterized and studied in great detail. The main references for discrete orthogonal polynomials are the books \cite{Niko} and \cite{Hyp}. See also \cite{AlvGMar}, \cite{DomMar1}, \cite{DomMar2}, and  \cite{SBFArea}. 

In this paper, we present a construction that produces for each polynomial family in  the q-Askey scheme   a  weight function  $w(x_k)$ that determines  a discrete orthogonality with respect to a quasi-definite moment functional.
The weights $w(x_k)$ are the values at $t=1$ of a sequence of basic hypergeometric functions $f_k(t)$ that are of type $_3\Phi_2$, or type $_2\Phi_1$. The functions $f_k(t)$ depend on the parameters that determine the coefficients in the  three-term recurrence relation of the corresponding polynomial family.  The basic hypergeometric functions $f_k(t)$ are convergent at $t=1$ and  we show how  they are obtained from  $f_0(t)$.

For some families  in the $q$-Askey scheme we obtain two discrete orthogonality weights associated with different sequences of nodes $x_k$.
In this paper we do not deal with the problem of characterizing the cases for which the weights are positive.

We obtain our results using the linear algebraic approach of our previous papers \cite{Mops}, \cite{Rec}, \cite{PSGH}, and \cite{Uni}.
In  \cite{Uni} we presented a unified construction of all the  hypergeometric and basic hypergeometric  orthogonal polynomial sequences that uses three linearly recurrent sequences  of numbers that satisfy certain difference equation of order three. The initial terms of such  sequences determine  the sequences of orthogonal polynomials and provide us with a uniform parametrization of all the  hypergeometric and basic hypergeometric sequences. In the present paper we use several results from \cite{Uni}.

In Section 2 we present some preliminary material related with  the construction of the hypergeometric $q$-orthogonal polynomials and their uniform parametrization from \cite{Uni}. In Section 3 we present  the construction of the weight functions obtained by  solving an infinite system of linear equations and use some  properties of basic  hypergeometric functions. In Section 4 we find general formulas for  the weight functions for the cases when some of the parameters are equal to zero. In Section 5 we present explicit expressions for  the weight functions for some families of polynomials in the $q$-Askey scheme.

\section{The class $\cH_q$ of basic  hypergeometric orthogonal polynomial sequences}

In this section we present some properties of  the class $\cH_q$ of the basic hypergeometric orthogonal polynomial sequences that were obtained in \cite{Uni}. 

Consider the homogeneous  difference equation
\begin{equation}\label{eq:diffeq}
 s_{k+3} = z ( s_{k+2} - s_{k+1}) + s_k,  \qquad k\ge 0,
\end{equation}
where $z=1+q + q^{-1}$.
The roots of the  characteristic polynomial  of  this equation  are  $1, q$, and $q^{-1}$.
In the present paper we consider only  the case with distinct roots. Therefore the general solution of \eqref{eq:diffeq}  is a  linear combination of the sequences $1, q^k$  and $q^{-k}$.  	Let the sequences  $x_k$, $h_k$, and $d_k$ be solutions of \eqref{eq:diffeq}. Then we can write these sequences as follows   
\begin{equation}\label{eq:xhe}
	x_k= b_0 + b_1 q^k  + b_2 q^{-k} , \qquad h_k = a_0 + a_1 q^k + a_2 q^{-k},\qquad  
	d_k= s_0+ s_1 q^k + s_2 q^{-k}. 
\end{equation}

The sequence $x_k$ determines the Newtonian basis $\{v_n(t): n\ge 0\}$ of the complex vector space of polynomials in the complex variable  $t$, defined by 
\begin{equation}\label{eq:Newton}
	v_n(t)= (t-x_0)(t-x_1)\cdots (t-x_{n-1}), \qquad n \ge 1,
\end{equation}
and $v_0(t)=1.$
We define the sequence $g_k$ by
\begin{equation} \label{eq:g}
g_k= x_{k-1} (h_k- h_0) + d_k, \qquad k \ge 1, 
\end{equation}
and $g_0=0$. Therefore we must have $d_0=0$ and hence $ s_0=-s_1-s_2$. In addition, we  suppose that $h_k \ne h_j$ whenever $k \ne j$,  and that  $g_k \ne 0$ for $k \ge 1$. 
We can take $a_0=0$ because $a_0$ does not appear in any of the formulas that we will obtain in our development.

Let $\cD$ be the linear operator on the space of polynomials defined  by
\begin{equation}\label{eq:operD}
\cD v_k = h_k v_k + g_k v_{k-1}, \qquad k \ge 0.
\end{equation}
	Since $g_0=0 $ we see that $\cD t^n$ is equal to $ h_n t^n $ plus a  polynomial of  degree less than $n$. For $n \ge 0$ let $u_n$ be a monic polynomial of degree $n$ which is an eigenfunction of $\cD$ with eigenvalue $h_n$. That is 
\begin{equation}\label{eq:eigenEq}
\cD u_k = h_k u_k, \qquad k \ge 0.
\end{equation}

The operator $\cD$ is a generalized difference operator which, in concrete examples, becomes a second order  difference  operator on a linear or quadratic $q$-lattice. In \cite{Uni} we showed that
\begin{equation}\label{eq:u}
 u_n(t) = \sum_{k=0}^n c_{n,k} v_k(t), \qquad n \ge 0,
\end{equation}
 where the coefficients $ c_{n,k}$ are given by
 \begin{equation}\label{eq:cnk}
c_{n,k}=\prod_{j=k}^{n-1} \frac{g_{j+1}}{h_n -h_j}, \qquad 0 \le k \le n-1,
 \end{equation}
and $c_{n,n}=1$ for $n \ge 0$. This expression for $u_n(t)$ was also  obtained by Vinet and Zhedanov in \cite{VZ} using a different approach. The idea of representing orthogonal polynomials in terms of a Newtonian basis was introduced by Geronimus in \cite{Ger}. 

The matrix of coefficients $C=[c_{n,k}]$ is an infinite lower triangular matrix where the coefficients of $u_n$ appear in the $n$-th row of $C$. Since $c_{n,n}=1$ for $n \ge 0$, the infinite matrix  $C$ is invertible. 
	Let $C^{-1}=[\hc_{n,k}]$ and define the polynomials  
\begin{equation}\label{eq:Newthk}
w_{n,k}(t)= \prod_{j=k}^{n-1} (t-h_j), \qquad 0 \le k \le n.
\end{equation}
Using basic properties of  divided differences we obtain	
\begin{equation} \label{Cinverse}
	 \hc_{n,k}=\frac{ \prod_{j=k+1}^n g_j}{ w_{n+1,k}\pr (h_k) }=\prod_{j=k+1}^n \frac{g_j}{h_k - h_j}, \qquad 0 \le k \le n-1, 
\end{equation}
	 and $\hc_{n,n}=1$ for $n \ge 0.$

The entries in the $0$-th column of $C^{-1}$ are given by
\begin{equation}\label{eq:GenMom}
 \hc_{k,0} =\prod_{j=1}^k \frac{g_j}{h_0-h_j}, \qquad k \ge 1, 
\end{equation}
and $\hc_{0,0}=1$. We denote them by $m_k=\hc_{k,0}$ for $ k \ge 0.$  They satisfy $m_0=1$ and
\begin{equation} \label{moments}
	\sum_{k=0}^n c_{n,k} \,  m_k = 0, \qquad n \ge 1.
\end{equation}
Note that the sequence $m_n$ satisfies a recurrence relation of order one. We will see that the numbers $m_k$ are the  generalized moments, with respect to the Newtonian basis $\{v_k(t) \}$, of a moment functional for which  the polynomial sequence $u_n(t)$ is orthogonal.

In \cite{Uni} we also proved that the polynomial sequence $u_n(t)$ satisfies a three-term recurrence relation of the form
\begin{equation}\label{eq:3term}
	u_{n+1}(t)= (t-\beta_n) u_n(t) - \alpha_n u_{n-1}(t), \qquad n \ge 1, 
\end{equation}
where the coefficients are given by
and 
\begin{equation}\label{eq:alpha}
 \alpha_n = \frac{g_n}{h_{n-1} -h_n} \left(\frac{g_{n-1}}{h_{n-2} - h_n} - \frac{g_n}{h_{n-1} - h_n} + \frac{g_{n+1}}{h_{n-1}-h_{n+1}} +x_n - x_{n-1} \right), 
\end{equation}
and 
\begin{equation}\label{eq:beta}
	\beta_n= x_n + \frac{g_{n+1}}{h_n - h_{n+1}} -\frac{g_n}{h_{n-1} - h_n}. \end{equation}

Writing $\alpha_n$ and $\beta_n$ in terms of the parameters introduced in \eqref{eq:xhe} we obtain
\begin{equation} \label{eq:alphabs}
	\alpha_n= \dfrac{(x-1) (a_1 x - q a_2) p_1(x) p_2(x)}{(a_1 x^2 -a_2) (a_1 x^2-q a_2)^2 ( a_1 x^2 - q^2 a_2)},\qquad n \ge 1,
\end{equation}
where $x=q^n$, 
$$p_1(x)= (a_1 x - a_2) (b_1 x^2 + b_0 q x + b_2 q^2) + q x ( s_1 x -s_2),$$
and 
$$p_2(x)= -\dfrac{a_1^2 x^3}{q^2 a_2} \, p_1\left(\dfrac{q a_2}{a_1 z}\right),$$
and
\begin{equation}\label{betabs}
	\beta_n= \dfrac{q b_0 (x^2-1) ( a_1^2 x^2- a_2^2) + d\,  x ( x-1) (a_1 x - a_2) - ( q s_1 -s_2) ( a_1 - q a_2) x^2}{(a_1 q x^2 -a_2) ( a_1 x^2 - q a_2)}, \quad n\ge 0,
\end{equation}
where
$$ d= (q+1) ( q a_1 b_2 + a_2 b_1) - q (s_1 + s_2 + ( a_1 + a_2) b_0).$$

Let us note that $\alpha_n$ and $\beta_n$  are rational functions of $q^n$ and are determined by the parameters $ a_1,a_2,b_0, b_1,b_2,s_1,s_2$. Note also that at least one of $a_1$ and $a_2$ must be different from zero.

It is easy to verify that $\alpha_n$ is invariant under the substitution $x\rightarrow {q a_2}/({a_1 x})$ and $\beta_n$ is invariant under the substitution $x \rightarrow {a_2}/({a_1 x})$.

The effect of the  exchange $q \leftrightarrow q^{-1}$ on the sequences $\alpha_n$ and $\beta_n$ is equivalent to the exchanges $a_1  \leftrightarrow a_2$,  $b_1  \leftrightarrow b_2$,  $s_1  \leftrightarrow s_2$.

The map 
$$(a_1,a_2, b_0,b_1,b_2,s_1,s_2) \rightarrow (\{\alpha_n\}_{n\ge 1}, \{\beta_n\}_{n\ge 0})$$
is not injective. Let us consider next an example that  shows that all the polynomial sequences in the $q$-Askey scheme correspond to a vector of parameters where at least one of $b_1$ and $b_2$ is not zero. The only exception is the sequence of continuous $q$-Hermite polynomials.

Let  $P=(a_1,a_2, b_0,0,0,s_1,s_2)$. We define the vectors   $P_1 =(a_1,a_2, b_0,b_1,0,{\hat s}_1,{\hat s}_2)$ where 
$$b_1=-\dfrac{b_0 a_2 +s_2}{a_2}, \qquad {\hat s}_1=\dfrac{q s_1 - b_0 a_2 - s_2}{q}, \qquad {\hat s}_2=- b_0 a_2, $$ 
and $P_2=(a_1,a_2, b_0,0, b_2,{\tilde s}_1,{\tilde s}_2)$, where
$$b_2=-\dfrac{b_0 a_1 + s_1}{a_1}, \qquad  {\tilde s}_1=-b_0 a_1,\qquad {\tilde s}_2=s_2 - q  a_1 b_0 - q s_1.$$
A straightforward computation shows that if $a_1\ne 0$ and $a_2\ne 0$ then the three vectors  $P, P_1, P_2$ 
 yield the same pair of coefficient sequences $(\alpha_n, \beta_n)$, if $a_1=0$ and $a_2\ne 0$ then $P$ and $P_1$ produce the same pair of coefficient sequences, and if $ a_1\ne 0$ and $a_2=0$ then $ P$ and $P_2$ produce the same pair  $(\alpha_n, \beta_n)$. 

By Favard's theorem \cite[Thm. 4.4]{Chi}, \cite{MarAlv}, if all the $\alpha_n$ are positive and the $\beta_n$ are real then  $\{u_n\}$ is orthogonal with respect to a positive-definite moment functional, and if all the $\alpha_n$ are nonzero then $\{u_n\}$ is orthogonal with respect to a  quasi-definite moment functional.

\section{Generalized moments and discrete orthogonality}

Let us suppose that the parameters that determine the sequences $x_k, h_k, $ and $d_k$ are such that 
$ h_k \neq h_j $ whenever $k \ne j$, and the coefficients $\alpha_n$, given by \eqref{eq:alpha}, are nonzero for $n \ge 1$. Then, by Favard's theorem   \cite[Thm. 4.4]{Chi}, \cite{MarAlv},    there exists a unique quasi-definite moment functional $\tau$ on the space of polynomials such that the polynomial sequence $\{u_n(t): n \ge 0\}$ is orthogonal with respect to $\tau$, that is,
\begin{equation}\label{orthogonality}
	\tau(u_n(t) u_m(t))= K_n \delta_{n,m}, \qquad K_n \ne 0, \ \ n,m \in \N.
	\end{equation}

  Since $u_0(t)=1$  we have
	\begin{equation}\label{eq:genmom}
	\tau(u_n(t)) =\sum_{k=0}^n c_{n,k} \tau(v_k(t) )=0, \qquad n \ge 1,
	\end{equation}
	and by the uniqueness of the inverse of the matrix  $C$ we must have  $\tau(v_k(t))=\hc_{k,0}=m_k$  for $k \ge 0$.
This means that the numbers $m_k$ are the generalized moments of $\tau$ with respect to the Newtonian basis $\{v_k(t): k\ge 0\}$. Note that $\tau(1)=m_0=1.$

\begin{theorem}\label{theorem1}
	Suppose that  $ a_2$ and $b_2$ are nonzero and that  $\alpha_n\ne 0$ for $n\ge 1$. Then there exists a unique weight function $w$, defined on the nodes $x_k$,  such that for every polynomial $p(t)$ we have 
\begin{equation}\label{eq:discrTau}
	\tau(p(t))=\sum_{k=0}^\infty p(x_k) w(x_k).
\end{equation}
\end{theorem}
 
{\it Proof:} From \eqref{eq:xhe} we see that the hypothesis imply that $\{x_k\}$ and $\{h_k\}$ are sequences of pairwise distinct numbers, and that the polynomial sequence $\{u_n(t): n\ge 0\}$ is orthogonal with respect to the quasi-definite  moment functional $\tau$.

We will write $r_j=w(x_j)$ for $j \ge 0$. Since  $\tau(v_k(t))=m_k$  for $k \ge 0$, the numbers $r_j$ that we want to find must satisfy 
\begin{equation}\label{eq:infsyst}
	m_k=\tau(v_k(t))= \sum_{j=0}^\infty  v_k(x_j)\, r_j , \qquad k \ge 0. 
\end{equation}
This is an infinite system of linear equations where the $r_j$ are unknown. Let us denote by $P$ the matrix of coefficients. Then  $P=[ v_k(x_j)]$, where $j$ is the index for rows, $k$ is the index for columns and $j$ and $k$ are non-negative integers. From the definition of the polynomials $v_k(t)$ in \eqref{eq:Newton} we can see that $P$ is an infinite  lower triangular matrix and its entries in the main diagonal are 
$$ v_k(x_k)=(x_k-x_0) (x_k-x_1) \cdots (x_k- x_{k-1}),$$
and hence  they are nonzero. Therefore $P$ is invertible. 

Using basic properties of divided differences with respect to the nodes $x_j$ it is easy to show that $P^{-1}$ is the lower triangular  matrix whose $(j,k)$ entry is  $1/v_{j+1}^\prime(x_k)$.
Therefore, from the system \eqref{eq:infsyst} we obtain
\begin{equation}\label{eq:weights}
	r_k=\sum_{j=k}^\infty \frac{m_j}{v_{j+1}^\prime(x_k)}, \qquad k\ge 0.
\end{equation}
We will show  that all the series  in \eqref{eq:weights} are convergent (or terminating) basic hypergeometric series.

Let us define the power series
\begin{equation}\label{eq:fkt}
	f_k(t)= \sum_{j=k}^\infty \frac{m_j}{v_{j+1}^\prime(x_k)} t^j, \qquad k \ge 0.
\end{equation}
Therefore $r_k=f_k(1)$ for $k \ge 0$.

From \eqref{eq:GenMom} we obtain  
\begin{equation}\label{eq:mj}
	m_j =\prod_{i=1}^j \frac{g_i}{h_0-h_i}, \qquad j \ge 1, 
\end{equation}
and $m_0=1$.

The denominators in  \eqref{eq:fkt} are 
\begin{equation}\label{eq:vprime}
	v_{j+1}^\prime(x_k)= \prod_{i=0, i\ne k}^j (x_k-x_i), \qquad k \ge  j \ge 0.
\end{equation}
Therefore $f_k(t)$ can be expressed as
\begin{equation}\label{eq:fktprod}
	f_k(t)= \sum_{j=k}^\infty \left(   \prod_{i=1}^j  \frac{g_i}{h_0-h_i} \prod_{i=0, i\ne k}^j \frac{1}{x_k-x_i}\right)    t^j, \qquad k \ge 0,
\end{equation}
and
\begin{equation}\label{eq:f0t}
	f_0(t)= 1 + \sum_{j=1}^\infty \left( \prod_{i=1}^j \frac{g_i}{ (h_0 - h_i) ( x_0 - x_i)} \right)  t^j.
\end{equation}

The coefficient of $t^j$, for $j \ge k$,   in the series $f_k(t)$ depends on the parameters $ a_1, a_2, b_0, b_1, b_2, s_1, s_2$. Let us introduce new parameters $ y_1,y_2, y_3, z_1, z_2$ as follows
\begin{equation}\label{eq:yz}
	\begin{split}
		&	a_1= \dfrac{z_1 a_2}{q}, \qquad b_1=\dfrac{z_2 b_2}{q}, \qquad  y_3=\dfrac{ z_1 z_2}{q  y_1 y_2},  \\ 
		&	s_1= \dfrac{-a_2}{q} \left( \left(  y_1 y_2 + y_1 y_3 + y_2 y_3 -\dfrac{z_2}{q} \right) b_2 + b_0 z_1\right),\\
		&	s_2= -a_2(  ( y_1+y_2 +y_3 - z_1) b_2 +b_0).  
	\end{split}
\end{equation}
The new parameters $z_1$ and $z_2$ are well defined by the first two equations because  $a_2\ne 0$ and $b_2\ne 0$.
Note that $y_3$ depends on $ y_1,y_2, z_1, z_2$. We use  the redundant parameter  $y_3$  because it simplifies the equations in \eqref{eq:yz} and other formulas that we will obtain next.

Let us denote by $\rho(k,j)$ the coefficient of  $t^j$  in the series $f_k(t)$.
By substitution in $\rho(0,j)$  of the parameters $a_1, b_1, s_1, s_2$, using the expressions in \eqref{eq:yz}, we obtain  
\begin{equation}\label{eq:coeff0yz}
	\rho(0,j)= \prod_{i=1}^j \frac{g_i}{ (h_0 - h_i) ( x_0 - x_i)} =\dfrac{(y_1;q)_j (y_2;q)_j (y_3;q)_j \, q^j}{(q;q)_j (z_1;q)_j (z_2;q)_j}, \qquad j \ge 1,
\end{equation}
where $(t;q)_k= (1-t) (1- q t) \cdots ( 1- q^{k-1} t)$  for $t \in \C$ and $k \in \N$. Note the factor $q^j$ in the numerator. Note also that only the parameters $z_1,z_2, y_1,y_2,y_3$ appear in the right-hand side. Thus we have essentially four parameters, since $y_3$ can be written in terms of the other four parameters.

The previous equation shows that $f_0(t)$ is a $_3\Phi_2$  basic hypergeometric series. If $|q| < 1$ or $|q|>1$ then, by the ratio test  the series $f_0(t)$ converges at $t=1$,   since $z_1 z_2= q y_1 y_2 y_3$. See  \cite[p. 5]{Gasper}.

 Using the change of parameters \eqref{eq:yz} we obtain
\begin{equation}\label{eq:rhokj}
	\rho(k,j)= \mu(k) \dfrac{(y_1;q)_j (y_2;q)_j (y_3;q)_j \, q^j}{(q;q)_{j-k}  (z_1;q)_j (q^k z_2;q)_j}, \qquad j \ge k \ge 1,
\end{equation}
where
\begin{equation}\label{eq:muk}
	\mu(k)=(-1)^k \dfrac{q^{\binom{k}{2}}} {(q;q)_k} \dfrac{(1-q^{2k-1} z_2)}{(1-q^{k-1} z_2)}, \qquad k\ge 1.
\end{equation}

The ratio test shows that the series $f_k(t)$ converges at $t=1$ if $|q|<1$ or $|q|>1$. 
This shows that the numbers $r_k$  defined in \eqref{eq:weights} are well defined and satisfy \eqref{eq:infsyst}. 
This completes the proof. \eop

 The coefficients of $f_0(t)$ in terms of the original parameters look like the coefficients of a $_3\Phi_2$ function, that is, like the right-hand side of  \eqref{eq:coeff0yz}. 
The formulas in \eqref{eq:yz}  for the change of parameters are obtained by substituting $a_1=q^{-1} z_1 a_2$ and $ b_1=q^{-1} z_2 b_2$ in  $f_0(t)$   and then  solving for $s_1,s_2$ and $y_3$  the system of  equations obtained from \eqref{eq:coeff0yz} by taking several values of $j$. 

The coefficients $\rho(k,j)$  of the series $f_k(t)$  satisfy
\begin{equation}\label{eq:sumrho}
\sum_{k=0}^j \rho(k,j) =0, \qquad j\ge 1,
\end{equation}
and
\begin{equation}\label{eq:recrho}
	\rho(k+1,j)=\dfrac{-q^k (1-q^{j-k}) ( 1-q^{2k+1} z_2) (1- q^{k-1} z_2)}{(1-q^{k+1}) ( 1-q^{2 k-1} z_2) (1-q^{j+k} z_2)} \, \rho(k,j), \qquad j >k.
\end{equation}

Let us  define the sequence of functions
\begin{equation}\label{eq:psi}
	\psi_k(t) = \mu(k) (z_2;q)_k \sum_{j=k}^\infty \dfrac{(q^{j-k+1};q)_k \, t^j}{(q^j z_2;q)_k}, \qquad k\ge 1.
\end{equation}
A simple computation shows that $f_k(t)$ is equal to the Hadamard product of the series $f_0(t)$ and $\psi_k(t)$, for $k\ge 1$.

The functions $f_k(t)$, for $k \ge 0$,  are invariant under the transformation
$$(q,y_1,y_2,y_3, z_1, z_2)\rightarrow(1/q, 1/y_1,1/y_2,1/y_3,1/z_1,1/z_2).$$

The change of parameters \eqref{eq:yz} applied to the coefficients of the three-term recurrence relation \eqref{eq:3term}
gives us 
\begin{equation}\label{eq:alphay}
	\alpha_n=\dfrac{ b_2^2  (x-1) (z_1 x -q^2) }{ z_1 (z_1 x^2 -q) (z_1 x^2 -q^2)^2 (z_1 x^2 -q^3)} \prod_{j=1}^3 ((y_j x -q) (z_1 x - q y_j)),
\end{equation}
where $x= q^n$,  and 
\begin{equation}\label{eq:betay}
	\beta_n=b_0 +  \dfrac{ b_2  x  (\epsilon_2 x^2 + \epsilon_1 x + \epsilon_0)   }{z_1 ( z_1 x^2 -1) ( z_1 x^2 -q^2)},
\end{equation}
where  $x= q^n$,  and
\begin{eqnarray*}
	\epsilon_2=& z_1(z_1^2 +( y_1+y_2+ y_3) q z_1 + (y_1 y_2 + y_1 y_3 + y_2 y_3) z_1 + q y_1 y_2 y_3), \\
	\epsilon_1=& - z_1 (q+1) (( q + y_1 + y_2 + y_3) z_1 + (y_1 y_2 + y_1 y_3 + y_2 y_3) q +  y_1 y_2 y_3), \\
	\epsilon_0=& q( z_1^2 +( y_1+y_2+ y_3) q z_1 + (y_1 y_2 + y_1 y_3 + y_2 y_3) z_1 + q y_1 y_2 y_3). 
\end{eqnarray*}

The entries of the matrix $C$, which are  the coefficients of the polynomials $u_n(t)$ with respect to the Newtonian basis $\{v_k(t):k\in \N\}$, become 
\begin{equation}\label{eq:cnkyz}
	c_{n,k}=\dfrac{b_2^{n-k}}{q^{(n-k) k}} \, \dfrac{(q^{k+1};q)_{n-k} (q^{k} y_1;q)_{n-k} (q^{k} y_2;q)_{n-k} (q^{k} y_3;q)_{n-k}}{(q;q)_{n-k} ( z_1 q^{n+k-1};q)_{n-k}},\qquad n\ge k \ge 0,
\end{equation}
in terms of the parameters $y_1,y_2,y_3, z_1, b_2$.

We consider next the family of Askey-Wilson polynomials with parameters $a,b,c,d$, as they are presented in \cite[Ch.14]{Hyp}. By a simple comparison of our formulas \eqref{eq:alphay} and \eqref{eq:betay} for $\alpha_n$ and $\beta_n$ with the corresponding formulas for the Askey-Wilson polynomials given  in \cite[Eq. 14.1.5]{Hyp},   we obtain
\begin{equation}\label{eq:AWyz}
	y_1=a b, \quad y_2=a c, \quad y_3=a d, \quad z_1=a b c d, \quad z_2=q a^2, \quad b_0=0, \quad b_2=\frac{1}{2 a}.
\end{equation}

Therefore the sequence of weights for the discrete orthogonality of the Askey-Wilson polynomials is given by
\begin{equation}\label{eq:AWweights}
	r_k=w(x_k)= \mu(k) \sum_{j=k}^\infty \dfrac{(ab;q)_j ( a c;q)_j ( a d;q)_j q^j}{(q;q)_{j-k} (a b c d;q)_j (q^{k+1} a^2;q)_j}, \qquad k \ge 0,
\end{equation}
where
\begin{equation}\label{eq:AWmuk}
	\mu(k)=(-1)^k  \dfrac{q^{\binom{k}{2}}}{(q;q)_k} \dfrac{(1- q^{2 k} a^2)}{(1-q^k a^2)}.
\end{equation}

The corresponding sequences of nodes $x_k$ and eigenvalues $h_k$ for the Askey-Wilson polynomials are given by
\begin{equation}\label{eq:AWxkhk}
	x_k=\dfrac{1}{2} \left( a \, q^k + \dfrac{ q^{-k}}{a}  \right), \qquad h_k= a_2 ( a b c d\, q^{k-1} + q^{-k}).
\end{equation}

We will show in the next section that Theorem \ref{theorem1} can be extended to the cases when at least one of $a_1$ and $a_2$ is nonzero and at least one of $b_1$ and $b_2$ is nonzero. 

\section{The cases with some parameters equal to zero}

In this section we consider the  cases with some of the four parameters $a_1,a_2, b_1,b_2$  equal to zero. Since our construction of the sequence of weights requires that $h_j \ne h_k$ and $ x_j \ne x_k$ whenever $j\ne k$, then at least one of $a_1, a_2$ must be nonzero and at least one of $b_1,b_2$ must be nonzero. Such restrictions give us nine cases.  
There are four cases with exactly one of the parameters equal to zero, four cases with two parameters equal to zero, and the case with no parameters equal to zero, that was considered in the previous section. Let us recall that all the pairs of sequences $(\alpha_n,\beta_n)$ of recurrence coefficients correspond to a vector of parameters that has  at least one of $b_1$ and $b_2$ different from zero.
It is clear that the results in the previous section apply to the   case where $a_1,a_2, b_1,b_2$ are nonzero numbers. We show next that they  can also be applied  to other three cases,  just  using some simple substitutions.

\subsection{ Case 1. $a_1=0$.}

If $a_1=0$ and $a_2,b_1,b_2$ are nonzero then, from \eqref{eq:yz} we obtain $z_1=0$ and $y_3=0$. Therefore the change of parameters  \eqref{eq:yz} becomes
\begin{equation}\label{eq:yza1}
	\begin{split}
		&	a_1= 0, \qquad b_1=\dfrac{z_2 b_2}{q}, \qquad  y_3=0,  \\ 
		&	s_1= \dfrac{a_2 b_2}{q^2} ( z_2 - q y_1 y_2   ) ,\\
		&	s_2= -a_2(  ( y_1+y_2 ) b_2 +b_0),  
	\end{split}
\end{equation}
and the coefficients $\rho(k,j)$ are
\begin{equation}\label{eq:rhoa1}
	\rho(k,j)=\mu(k) \dfrac{(y_1;q)_j (y_2;q)_j q^j }{(q;q)_{j-k} (q^j z_2;q)_j}, \qquad j \ge k,
\end{equation}
where $\mu(k)$ is defined in \eqref{eq:muk}. 
Therefore $f_k(t)$ is a $_2\Phi_1$ basic hypergeometric function.

In this case the eigenvalues are  $h_k=a_2 q^{-k}$, the nodes are    $x_k=b_0 + b_2(q^{k-1} z_2 +q^{-k})$ and the coefficients in the three-term recurrence relation are  
\begin{equation}\label{eq:alpha1}
\alpha_n= b_2^2 q^{-3} (x-1) ( x y_1 y_2 -z_2) ( x y_1 -q) ( x y_2 -q),
\end{equation}
where $x=q^n$, and
\begin{equation}\label{eq:beta1}
	\beta_n=b_0  -b_2 q^{-1} x (y_1 y_2 (1+q) x -z_2 -q (y_1+y_2) -y_1 y_2).
\end{equation}

\subsection{ Case 2. $b_1=0$.}

If $b_1=0$ then $z_2=0$ and $y_3=0$  and the change of parameters becomes 
\begin{equation}\label{eq:yzb1}
	\begin{split}
		&	a_1= \dfrac{z_1 a_2}{q}, \qquad b_1=0, \qquad  y_3=0,  \\ 
		&	s_1=-a_2 q^{-1} (b_2 y_1 y_2 + z_1 b_0),\\
		&	s_2= a_2(( b_2 (z_1-y_1 -y_2) -b_0),  
	\end{split}
\end{equation}
and the coefficients $\rho(k,j)$ are
\begin{equation}\label{eq:rhob1}
	\rho(k,j)=  (-1)^k \dfrac{q^{\binom{k}{2}}}{(q;q)_k}  \dfrac{(y_1;q)_j (y_2;q)_j q^j }{(q;q)_{j-k} ( z_1;q)_j}, \qquad j \ge k.
\end{equation}
Therefore in this case  $f_k(t)$ is also  a $_2\Phi_1$ basic hypergeometric function.

The eigenvalues are  $h_k=a_2 ( z_1 q^{k-1} + q^{-k})$, the nodes are    $x_k=b_0 + b_2 q^{-k}$ and the coefficients in the three-term recurrence relation are  
\begin{equation}\label{eq:alphb1}
	\alpha_n= \dfrac{b_2^2 q x (x-1) (z_1 x - q^2) ( z_1 x - q y_1)  (z_1 x - q y_2) (y_1 x -q) ( y_2 x -q)}{(z_1 x^2 -q) ( z_1 x^2 - q^2)^2 ( z_1 x^2 - q^3)},
\end{equation}
where $x=q^n$, and
\begin{equation}\label{eq:betb1}
	\beta_n=b_0 + \dfrac{b_2 x (\epsilon_2 x^2 + \epsilon_1 x + \epsilon_0) }{ (z_1 x^2 -1) (z_1 x^2 - q^2)},
\end{equation}
where 
\begin{eqnarray*}
	\epsilon_2=&  z_1 (z_1 +( y_1+y_2) q  + y_1 y_2   ), \\
	\epsilon_1=&  - (q+1) (( q + y_1 + y_2 ) z_1 + (y_1 y_2 ) q ), \\
	\epsilon_0=&  q( z_1 +( y_1+y_2) q  + y_1 y_2   ). 
\end{eqnarray*}

\subsection{Case 3.  $a_1=0$  and $b_1=0$.}

If  $a_1=0$  and $b_1=0$ the  change of parameters becomes
\begin{equation}\label{eq:yzb11}
	\begin{split}
		&	a_1= 0, \qquad b_1=0, \qquad  y_3=0,  \\ 
		&	s_1=-a_2 b_2  q^{-1}  y_1 y_2,\\
		&	s_2=- a_2(( b_2 (y_1 +y_2) +b_0),  
	\end{split}
\end{equation}
and the coefficients $\rho(k,j)$ are
\begin{equation}\label{eq:rhoa1b1}
	\rho(k,j)= (-1)^k  \dfrac{q^{\binom{k}{2}}}{(q;q)_k} \dfrac{(y_1;q)_j (y_2;q)_j q^j}{(q;q)_{j-k}}, \qquad j\ge k.
\end{equation}

The eigenvalues are $h_k= a_2 q^{-k}$,  the nodes are $x_k= b_0 + b_2 q^{-k}$, and the recurrence coefficients are
\begin{equation}\label{eq:alpha1b1}
\alpha_n= b_2^2 q^{-3} y_1 y_2 x (x-1) (y_1 x - q) ( y_2 x - q),
\end{equation}
and
\begin{equation}\label{eq:beta1b1}
\beta_n=b_0 -b_2 q^{-1} x ( y_1 y_2 (q+1) x - q ( y_1+y_2) - y_1 y_2).
\end{equation}

\subsection{Case 4.  $a_2=0$.}

We substitute $a_2=0$ and $b_1=z_2 b_2/q $ in $f_0(t)$ and observe that the coefficient $\rho(0,j)$  of $t^j$ is of the form 
\begin{equation}\label{eq:rho4}
\rho(0,j)=\dfrac{(y_1;q)_j (y_2;q)_j r^j}{(q;q)_j (z_2;q)_j},
\end{equation}
for some numbers $y_1,y_2, p $. Then we solve for $s_1,s_2, y_2$ the system of equations obtained from \eqref{eq:rho4} by  taking $j=1,2,3. $  We obtain the change of parameters
\begin{equation}\label{eq:yza2}
	\begin{split}
		&	a_2= 0, \qquad b_1= \dfrac{ z_2 b_2}{q} , \qquad  y_2=\dfrac{z_2}{ p y_1},  \\ 
		&	s_1=-\dfrac{a_1}{q} \left(p y_1 +\dfrac{z_2}{y_1}\right)b_2 - a_1 b_0,   \\
		& 	s_2=  a_1  b_2 (q-p).  
	\end{split}
\end{equation}
This change of parameters gives us 
\begin{equation}\label{eq:rhoa2}
	\rho(k,j)= (-1)^k  \dfrac{q^{\binom{k}{2}}  (1- q^{2k-1} z_2)}{(q;q)_k (1- q^{k-1} z_2)} \dfrac{(y_1;q)_j \left(\dfrac{z_2}{p y_1};q \right)_j \, p^j }{(q;q)_{j-k}  (q^k z_2;q)_j}, \qquad j\ge k.
\end{equation}

In this case the recurrence coefficients are 
\begin{equation}\label{eq:alpha2}
	\alpha_n= \dfrac{b_2^2}{q y_1} \dfrac{(x-1) ( x-p)  (x y_1 -q) ( x z_2 - q p  y_1)}{x^4},
\end{equation}
and
\begin{equation}\label{eq:beta2}
	\beta_n=b_0 + \dfrac{b_2}{q y_1} \dfrac{((q+y_1) p + q) y_1 x - p y_2 (1+q) + z_2   }{x^2},
\end{equation}
where $x=q^n$,  the eigenvalues are $h_k= a_1 q^k$, and the nodes are $x_k= b_0 + b_2 (z_2 q^{k-1}  +q^k)$.

\subsection{Case 5.  $a_2=0$ and $b_1=0$.}

We put $b_1=0$ and $z_2=0$ in the formulas of the previous case and obtain the change of parameters
\begin{equation}\label{eq:yza2b1}
	\begin{split}
		&	a_2= 0, \qquad b_1= 0, \qquad z_2=0,  \\ 
		&	s_1=-\dfrac{a_1}{q} (p y_1 b_2  + q b_0),   \\
		&	s_2=  a_1  b_2 (q-p),  
	\end{split}
\end{equation}
which gives us
\begin{equation}\label{eq:rhoa2b1}
	\rho(k,j)= (-1)^k \dfrac{ q^{\binom{k}{2}} }{(q;q)_k} \dfrac{(y_1;q)_j p^j}{(q;q)_{j-k}}.
\end{equation}
In this case the recurrence coefficients become
\begin{equation}\label{eq:alpha2b1}
	\alpha_n= b_2^2 p  \dfrac{(x-1) (x-p)  (q- x y_1)}{x^4},
\end{equation}

\begin{equation}\label{eq:beta2b1}
	\beta_n=\frac{b_2}{q} \dfrac{(p ( q+y_1) +q) x - p (q+1) }{x^2},
\end{equation}
where, as usual, $x=q^n$.
The eigenvalues are $h_k= a_1 q^k$ and the nodes $x_k=b_0 + b_2 q^{-k}$.

\subsection{Case 6. $b_2=0$.}
We proceed as in Case 4, but now we put $b_2=0$ and $a_1= z_1 a_2/q$ in $f_0(t)$, and then we see that the coefficients $\rho(0,j)$ are of the form
\begin{equation}\label{eq:rho6}
	\rho(0,j)= \dfrac{(y_1;q)_j (y_2;q)_j \, p^j}{(q;q)_j (z_1;q)_j }.
\end{equation}
We solve for $ s_1,s_2, y_2$ the system of equations obtained by taking $j=1,2,3$ in \eqref{eq:rho6} and obtain
the change of parameters
\begin{equation}\label{eq:yzb2}
	\begin{split}
		&	a_1=\dfrac{ z_1 a_2}{q}, \qquad b_2 = 0, \qquad y_2= \dfrac{z_1}{p y_1}  \\ 
		&	s_1=-\dfrac{a_2}{q y_1}    \left( p b_1 y_1^2 +(b_0 z_1 - b_1) y_1 + z_1 b_1 \right),   \\
		&	s_2= -  a_2 (b_0 +p   b_1),  
	\end{split}
\end{equation}
which gives us
\begin{equation}\label{eq:rhob2}
	\rho(k,j)= (-1)^k \dfrac{ q^{\binom{k}{2}} }{(q;q)_k} \dfrac{(y_1;q)_j \left(\dfrac{z_1}{ p y_1} ; q \right)_j  \, p^j }{ q^{k(j-1)}  (q;q)_{j-k} (z_1;q)_j}.
\end{equation}
In this case the recurrence coefficients become
\begin{equation}\label{eq:alphb2}
	\alpha_n= - \dfrac{- q b_1^2}{y_1^2} \dfrac{  x (x-1) (x z_1-q^2) (x z_1 -q y_1) ( x p y_1-q) ( x y_1-q) ( x z_1 - q p y_1)     }{(x^2 z_1 -q) (x^2 z_1 - q^2)^2  (x^2 z_1 -q^3)  },
\end{equation}

\begin{equation}\label{eq:betab2}
	\beta_n=b_0- \frac{b_1}{y_1} \dfrac{x ( z_1 \epsilon_1 x^2 - ( q+1) \epsilon_2 x + q \epsilon_1     }{(x^2 z_1-1) (x^2 z_1-q^2)},
\end{equation}
where, as usual, $x=q^n$ and 
$$\epsilon_1= p y_1^2 + q  y_1 ( p+1) + z_1 y_1, \qquad \epsilon_2= q p y_1^2 + z_1 y_1 ( p+1) + q z_1. $$

The eigenvalues are $h_k= a_2 (z_1 q^{k-1} + q^{-k})$ and the nodes $x_k=b_0 + b_1 q^{k}$.

\subsection{Case 7.  $ b_2=0$  and $a_1=0$.}

The formulas for this case are obtained from the corresponding formulas in Case 6, just putting $a_1=0$ and $z_1=0$.

The change of parameters becomes
\begin{equation}\label{eq:yzb2a1}
	\begin{split}
		&	a_1= 0, \qquad b_2 = 0, \qquad z_1= 0  \\ 
		&	s_1=-\dfrac{a_2 b_1}{q}  (p y_1-1 ),   \\
		&	s_2= -  a_2 (b_0 +p  b_1),  
	\end{split}
\end{equation}
which gives us
\begin{equation}\label{eq:rhob2a1}
	\rho(k,j)= (-1)^{k}  \dfrac{ q^{\binom{k}{2}} }{  (q;q)_k} \dfrac{(y_1;q)_j  p^j }{ q^{k(j-1)}  (q;q)_{j-k} }.
\end{equation}
In this case the recurrence coefficients become
\begin{equation}\label{eq:alphb2a1}
	\alpha_n=  \dfrac{p  b_1^2}{q^3}  x (x-1) (x p  y_1 -q) (x y_1 - q),
\end{equation}

\begin{equation}\label{eq:betab2a1}
	\beta_n=b_0- \frac{b_1}{q}  x ((q+1) p y_1 x -q ( p+1) - p y_1 ), 
\end{equation}
where, as usual, $x=q^n$.
 The eigenvalues are $h_k= a_2  q^{-k}$ and the nodes $x_k=b_0 + b_1 q^{k}$.

\subsection{Case 8. $a_2=0$ and $b_2=0$.}

If $a_2=0$ and $b_2=0$ it is easy to see that the coefficients $\rho(0,j)$  of the series $f_0(t)$ are of the form
\begin{equation}\label{eq:rho8}
	\rho(0,j) = (-1)^j q^{-\binom{j}{2}} \dfrac{(y_1;q)_j (y_2;q)_j p^j}{(q;q)_j}.
\end{equation}
Solving for $s_1,s_2, y_2$ the system of equations obtained from \eqref{eq:rho8} by  taking $j=1,2,3$ we obtain the change of parameters
\begin{equation}\label{yza2b2}
	\begin{split}
		&	a_2= 0, \qquad b_2 = 0, \qquad y_2=\dfrac{1}{p y_1},  \\ 
		&	s_1=-\dfrac{a_1 (b_1 p   y_1^2 + b_0 y_1 +b_1)}{y_1},   \\
		&	s_2= - q  p a_1   b_1.  
	\end{split}
\end{equation}
This change of parameters gives us
\begin{equation}\label{rho8}
	\rho(k,j) = (-1)^{k+j} \dfrac{q^{\binom{k}{2}}}{(q;q)_k} \dfrac{(y_1;q)_j \left(\dfrac{1}{p y_1};q \right)_j \, p^j }{q^{\binom{j}{2}} q^{k(j-1)} (q;q)_{j-k}}.
\end{equation}
The recurrence coefficients are 
\begin{equation}\label{eq:alpha2b2}
	\alpha_n= -\dfrac{q p  b_1^2 (x-1) ( x- q p  y_1) ( x y_1-q)}{ y_1 x^4},
\end{equation}
and
\begin{equation}\label{eq:beta2b2}
\beta_n=b_0 + \dfrac{ b_1 (p y_1^2 + q p  y_1 +1) x - (q+1) p  y_1 }{y_1 x^2},
\end{equation}
where $x=q^n$.

The eigenvalues are $h_k= a_1 q^k$ and the nodes are $x_k= b_0 + b_1 q^k$.

\section{Some examples.} 

Since the   nodes $x_k$ are defined  by $x_k= b_0 + b_1 q^k + b_2 q^{-k}$, for $k\ge 0$,    
from \eqref{eq:Newton} it is easy to see that a change in  the constant term $b_0$  corresponds to a translation of the polynomials $v_n(t)$ in the Newton-type basis,  and therefore to  the same translation of the orthogonal polynomials $u_n(t)$. We will consider that the polynomial sequences $u_n(t)$ and ${\tilde u}_n(t)$ are equivalent if there exists a constant $\lambda$ such that ${\tilde u}_n(t)=u_n(t+ \lambda)$, for $n \ge 0.$
 When the coefficients  $\beta_n$ of the three-term recurrence relation are expressed in terms of the parameters $ y_1, y_2, z_1,z_2$  the parameter $b_0$ appears as an additive constant. 

Let us recall that 
$$f_k(t)= \sum_{j=k}^\infty \rho(k,j) \, t^j, \qquad k\ge 0, $$
and that the weights associated with the nodes $x_k$  are the numbers $f_k(1)$, for $k \ge 0$.

We present next   explicit formulas for the coefficients $\rho(k,j)$  for some families of orthogonal polynomials in the $q$-Askey scheme.

\subsection{Continuous big $q$-Hermite,  \cite[14.18]{Hyp}.}
For this family we have 
$$\alpha_n=\dfrac{1-q^n}{4}, \qquad \beta_n=\dfrac{a q^n}{2}, \qquad x_k=\dfrac{a q^k + (a q^k)^{-1}}{2},$$
and 
$$z_1=0, \qquad z_2= a^2 q, \qquad y_1=0,\qquad  y_2=0, \qquad y_3=0.$$
These values of the parameters give us
$$\rho(k,j)= \dfrac{(-1)^k q^{\binom{k}{2}} (1- q^{2 k} a^2) \, q^j}{(1-q^k a) (q;q)_k (q;q)_{j-k} (q^{k+1} a^2;q)_j},\qquad j \ge k\ge 0.$$

\subsection{Discrete $q$-Hermite I, \cite[14.28]{Hyp}.}
In this case we have
$$\alpha_n=\dfrac{q^n ( 1-q^n)}{q}, \quad \beta_n=0, \quad x_k=q^k,$$
and the parameters are 
$$a_1=0,\quad b_0=0, \quad  b_1=1, \quad b_2=0, \quad p=-1, \quad y_1=0, \quad s_1=\dfrac{a_2}{q}, \quad s_2=a_2, $$
where $a_2$ is an arbitrary nonzero number.  Since $a_1=0$ and $b_2=0$ this example is in Case 7 and we obtain
$$\rho(k,j)=\dfrac{(-1)^{k+j} q^{\binom{k}{2}} }{q^{k(j-1)} (q;q)_k (q;q)_{j-k}}, \qquad j \ge k \ge 0.$$

\subsection{Little $q$-Jacobi, \cite[14.12]{Hyp}.}

The normalized three-term  recurrence relation for the little $q$-Jacobi polynomials is given in \cite[14.12.4]{Hyp}, with parameters $a$ and $b$.
The recurrence coefficients can be obtained in  three different ways in terms of the parameters $b_0, b_1,b_2,a_1,a_2,s_1,s_2$.
The first one is
$$a_1= a b q a_2, \quad b_0=0, \quad b_1=0, \quad b_2=\dfrac{1}{ q b}, \quad s_1= 0, \quad s_2=(q a -1) a_2,$$
where $a_2$ is an arbitrary nonzero number. Since $b_1=0$ and $a_1, a_2, b_2$ are nonzero this is an example of Case 2 of the previous section. The sequence of nodes is $x_k= b^{-1} q^{-k-1}$ and the coefficients $\rho(k,j)$ of the weights
are given by \eqref{eq:rhob1}  taking $y_1=b q, y_2=0, z_1=q^2 a b, z_2=0$. We obtain
$$\rho(k,j)=\dfrac{(-1)^k q^{\binom{k}{2}} (q b; q)_j q^j}{(q;q)k (q;q)_{j-k} (q^2 a b;q)j}, \qquad j \ge k. $$

The second set of values of the  parameters is
$$ a_1= a b q a_2, \quad b_0=0, \quad b_1=1, \quad b_2=0, \quad s_1= - q^{-1}(q a -1) a_2, \quad  s_2=0,$$
where $a_2$ is an arbitrary nonzero number. Since $b_2=0$ and $a_1, a_2, b_1$ are nonzero this is an example of Case 6  of the previous section. The sequence of nodes is $x_k= q^k$ and the coefficients $\rho(k,j)$ of the weights
are given by \eqref{eq:rhob2}  taking $y_1= b q, z_1=q^2 a b, y_2=0, p=0 $. We obtain 
$$\rho(k,j)= (-1)^{k+j} \dfrac{q^{\binom{k}{2}} q^{\binom{j}{2}} (q b ;q)_j a^j q^j}{q^{k(j-1)} (q;q)_j (q;q)_{j-k} (q^2 a b;q)_j}, \qquad j \ge k \ge 0.$$

The third set of values of the parameters is 
$$a_1=-q a b, \quad a_2=-1, \quad b_0=0,\quad b_1=0, \quad b_2=0, \quad s_1=a, \quad s_2=1.$$
Since the nodes are $x_k=0$, for $k \ge 0$, we do not obtain a discrete orthogonality with this set of values of the parameters.

\subsection{ $q$-Laguerre, \cite[14.21]{Hyp}.}

This family is obtained taking 
$$a_2=0, \quad b_0=1, \quad b_1=0, \quad b_2=-1, \quad s_1=0, \quad s_2=a_1 (q^{-{\alpha}} -q), $$
where $a_1$ is an arbitrary nonzero number.

The coefficients $\rho(k,j)$ are in this case given by
$$\rho(k,j)=\dfrac{(-1)^k q^{\binom{k}{2}} q^{-j \alpha}   }{(q;q)_k (q;q)_{j-k}    }, \qquad j \ge k \ge 0.$$
The sequence of nodes is $x_k=-q^{-k}$, for $k \ge 0$. Note that $\rho(0,j)$ are the coefficients of the $q$-exponential function of $q^{-\alpha}$. 

\subsection{Dual $q$-Hahn, \cite[14.7]{Hyp}.}

The dual $q$-Hahn polynomials can be obtained with three different sets of values of the parameters. We consider first the case with
$$a_1=0, \quad b_0=0, \quad b_1=q^{-N}, \quad b_2=q \gamma \delta q^N, \quad s_1= a_2 q^{-N} ( q^{-1} -\gamma), \quad s_2=- q a_2 \gamma (1+ \delta),$$
where $a_2$ is an arbitrary nonzero number and $N, \gamma, \delta$ are the parameters used in  \cite[14.7]{Hyp} to describe the $Q$-Hahn polynomials. This example is in Case 1 of the previous section. Therefore the coefficients $\rho(k,j)$ can be obtained from \eqref{eq:rhoa1} taking 
$$y_1= q^{-N},\qquad  y_2=\dfrac{q^{-N}}{\delta} ,\qquad  z_2=\dfrac{q^{-2 N}}{\gamma \delta},$$
which are the values of $y_1,y_2,z_2$ obtained from \eqref{eq:yza1} with the values of the parameters given above.
Since $y_1=q^{-N}$ we obtain a discrete orthogonality on the nodes $x_k= \gamma \delta q^{N+1-k} + q^{k-N}$, for $0 \le k \le N$.

The second set of values of the parameters is
$$a_1=0, \quad b_0=0, \quad b_1=q \gamma, \quad b_2= \delta, \quad s_1=a_2 \gamma ( 1-q^{-N}), \quad s_2=- a_2 (\gamma \delta q + q^{-N}),$$
where $a_2$ is an arbitrary nonzero number. The nodes are $x_k= \gamma q^{k+1} + \delta q^{-k}$, for $k \ge 0$.  We obtain the coefficients $\rho(k,j)$ from \eqref{eq:rhoa1} taking 
$$y_1= \gamma q, \qquad y_2=\dfrac{q^{-N}}{\delta}, \qquad z_2= \dfrac{\gamma q^2}{\delta}.$$

The third set of values of the parameters is 
$$a_1=0,\quad b_0=0, \quad b_1=\gamma \delta q, \quad b_2=1, \quad s_1=a_2 \gamma (\delta - q^{-N}), \quad s_2= -a_2 ( \gamma q + q^{-N}).$$
Substitution of these values in \eqref{eq:yza1} yields
$$y_1= q^{-N}, \qquad y_2= \dfrac{q^{-N}}{\delta}, \qquad z_2= \dfrac{q^{-2 N}}{\gamma \delta},$$
and then \eqref{eq:rhoa1} gives us the corresponding coefficients $\rho(k,j)$. The discrete orthogonality obtained in this case is essentially  equivalent to the one described in \cite[14.7.2]{Hyp}. In this case the nodes are $x_k=\gamma \delta q^{k+1} + q^{-k}$, for $k \ge 0$.

\section{Final remarks.}
We have shown that all the families of basic hypergeometric polynomial sequences in the $q$-Askey scheme have al least one discrete orthogonality on a sequence of nodes of the form $x_k=b_0 + b_1 q^k + b_2 q^{-k}$, where at least one of $b_1$ and $b_2$ is nonzero. The only exception is the sequence  of continuous $q$-Hermite polynomials. Some families are orthogonal on two different sequences of nodes $x_k$.

In some cases the basic hypergeometric series $f_k(t)$ are obtained from $f_0(t)$ by repeated application of the  $q$ derivative or the  $q^{-1}$ derivative. In some of the simplest cases the functions $f_k(t)$ are related with the $q$-exponential function.
 We have not studied the problem of characterizing the cases that have positive weights 

The formulas \eqref{eq:alphay} and \eqref{eq:betay} for the recurrence coefficients $\alpha_k$ and $\beta_k$ in terms of the parameters $y_1,y_2, z_1, z_2$ can be used to construct a classification of the $q$-orthogonal polynomial sequences that would improve the $q$-Askey scheme.

\end{document}